\documentclass[12pt]{article}
\usepackage{natbib}
\usepackage{dcolumn}
\usepackage{graphicx}
\usepackage{amsmath}
\usepackage{amsfonts}
\usepackage{amsopn}
\usepackage{amsthm}
\usepackage[figuresright]{rotating}
\usepackage{amsmath,amsthm,amsfonts,epsfig,amssymb,natbib,eucal,eufrak}
\usepackage{enumitem}
\usepackage{color}
\newtheorem{thm}{Theorem}

\newtheorem{cor}{Corollary}

\begin{document}
\begin{center}
{\Large\bf On the structure of UMVUEs}
\end{center}

\vspace{.5cm}
\centerline{Abram M. Kagan${}^{a}$ and Yaakov Malinovsky${}^{b, *}$}
\begin{center}
${}^{a}$ Department of Mathematics, University of Maryland, College Park, MD 20742, USA \\
${}^{b}$ Department of Mathematics and Statistics , University of Maryland, Baltimore County, Baltimore, MD 21250, USA
\end{center}
\centerline{${}^{a}$   {\it email}: amk@math.umd.edu}
\centerline{${}^{b}$  {\it email}: yaakovm@umbc.edu, ${}^{*}$corresponding author}

\vspace{.5cm} \begin{abstract}
In all setups when the structure of UMVUEs is known, there exists a subalgebra $\cal U$ (MVE-algebra) of the basic $\sigma$-algebra such that all  $\cal U$-measurable statistics with finite second moments are UMVUEs. It is shown that MVE-algebras are, in a sense, similar to the subalgebras generated by complete sufficient statistics. Examples are given when these subalgebras differ, in these cases a new statistical structure arises.
\end{abstract}
\vspace{9pt} \noindent {\it Key words:} Categorical Data; Completeness; MVE-algebra; Sufficiency.\\
\vspace{9pt} \noindent {\it 2000 Mathematics Subject Classification}. Primary 62B99; Secondary 62F10, 62G05.

\section{Introduction}
Among C.~R. Rao's major contributions to the foundations of statistical inference (Cram$\acute{e}$r-Rao inequality, Rao-Blackwellization, Rao's score test, to name the best known), the following observation in \cite{R1952} is not widely known.

Let $\displaystyle \left(\mathcal{X}, \mathcal{A}, \mathcal{P}\right)$ be a standard statistical model with $\displaystyle \left(\mathcal{X}, \mathcal{A}\right)$ a measurable space, and let $\mathcal{P}=\{P_{\theta}, \theta \in \Theta\}$ be a family of probability distributions of a random element $\displaystyle X \in \left(\cal{X}, \cal{A}\right)$
parameterized by a general parameter $\theta$.

A statistic $\widehat{g}(X)$ with $E_{\theta}|\widehat{g}(X)|^2<\infty,\,\theta \in \Theta$ is called a UMVUE (uniformly minimum variance unbiased estimator), more precisely, the UMVUE of $g(\theta)=E_{\theta}\widehat{g}(X)$ if
$$Var_{\theta}\left(\widehat{g}(X)\right)\leq Var_{\theta}\left(\widetilde{g}(X)\right),\,\,\forall\, \theta \in \Theta$$
for any statistic $\widetilde{g}(X)$ (with finite second moment) with $E_{\theta}\widetilde{g}(X)=E_{\theta}\widehat{g}(X),\,\,\theta \in \Theta$.

If $\displaystyle \widehat{g}_1\left(X\right)$ and $\displaystyle \widehat{g}_2\left(X\right)$ are UMVUE's, so is their linear combination $\displaystyle c_1\widehat{g}_1\left(X\right)+c_2\widehat{g}_2\left(X\right),\,\,\,c_1, c_2 \in \mathbb{R}$. Rao observed that their product  $\displaystyle \widehat{g}\left(X\right)=\widehat{g}_1\left(X\right)\widehat{g}_2\left(X\right)$ is also a UMVUE provided that
$\displaystyle E_{\theta}|\widehat{g}\left(X\right)|^2<\infty,\,\,\,\theta\in \Theta$ (this is not guaranteed by
$\displaystyle E_{\theta}|\widehat{g}_{1}\left(X\right)|^2+E_{\theta}|\widehat{g}_{2}\left(X\right)|^2<\infty,\,\,\,\theta\in \Theta$).
In all known (to the authors) setups{\color{blue},} the class of UMVUE's has the following structure: there exists a subalgebra (an MVE-algebra) $\displaystyle \widehat{\cal A} \subseteq {\cal A}$ such that a statistic $\displaystyle \widehat{g}\left(X\right)$ is a UMVUE
if and only if $\displaystyle \widehat{g}\left(X\right)$ is $\cal \widehat{A}$-measurable.
Combining all this leads to a hypothesis that such is the structure of UMVUEs in general. The results in the paper are closely related to this hypothesis, which, if true, introduces a new statistical structure.

In Section \ref{sec:2} a few setups with known structure of UMVUEs are presented. The first one (see Subsection \ref{subsec:21}) is classical, in a sense, of a family $\displaystyle \cal P$ possessing a complete sufficient statistic/subalgebra.
The MVE-algebra is complete sufficient subalgebra and every estimable parametric function possesses a UMVUE. The other setups, when the minimal sufficient subalgebra is incomplete, are less known. Plainly, due to the Rao-Blackwell theorem, MVE-algebras are subalgebras of the minimal sufficient subalgebras. In Subsection \ref{subsec:22} the setup of a partial complete sufficient subalgebra is considered, that formally generalizes the setup of Subsection \ref{subsec:21} (for a relation between partial completeness and completeness see a recent paper  \cite{KMM2014}). In Subsection \ref{subsec:23} a geometric construction leads to a description of the MVE-algebras in case of categorical $X$. Lehmann's example presented in Subsection \ref{subsec:24} deals with discrete $X$ taking countable many values and thus is not covered by the results of Subsection \ref{subsec:23}.

In Section \ref{sec:3} some results for arbitrary (not necessarily generated by complete sufficient statistics) MVE-algebras are proved{\color{blue},} showing that the MVE-algebras preserve properties of complete sufficiency.

The existence of UMVUEs in the Nile problem by R.~Fisher (see \cite{F1936}, \cite{F1973}) and related problems attracted recently some attention (see \cite{NS2012}). For its partial solution see \cite{KM2013}.
For the relation between sufficient statistic, and sufficient subalgebras we refer to \cite{B1957}.

\section{Setups with known structure of UMVUEs}
\label{sec:2}
A statistic $\chi(X)$ is called an unbiased estimator of zero (or zero mean statistic) if
$$E_{\theta}\chi(X)=0,\,\,\,\forall\, \theta \in \Theta.$$
A well known result (see, e.g., \cite{Lehmann1998}, p. 85) gives a necessary and sufficient conditions for a statistic $\widehat{g}(X)$ to be a UMVUE.
%\begin{lem}
%\label{le:1}
A statistic $\widehat{g}(X)$ is a UMVUE if and only if it is uncorrelated with any zero mean statistic $\chi(X)$ with $E_{\theta}|\chi(X)|^{2}<\infty,\,\,\forall\,\, \theta \in \Theta$.
%\end{lem}
Its immediate corollary is that the class of UMVUEs is a linear space: if $\widehat{g}_1$ and $\widehat{g}_2$ are UMVUEs, so is $\widehat{g}(X)=c_1\widehat{g}_1(X)+c_2\widehat{g}_2(X)$ for any constants $c_1,\,c_2$.

In the following subsections, setups with known structure of the class of UMVUEs are presented.
\subsection{Families with complete sufficient subalgebras/statistics}
\label{subsec:21}
Let a statistic $\displaystyle S=S(X),\,S: \left(\cal{X}, \cal{A}\right)\rightarrow \left(\cal{S}, \mathcal{B}\right)$
be a complete sufficient statistic for $\cal P$, i.e., $\displaystyle  E_{\theta}\chi\left(S(X)\right)=0,\,\forall\, \theta \in \Theta$ implies
\begin{align}
\label{eq:com}
P_{\theta}\left(\chi=0\right)=1\,\,\,\,\text{for all}\,\,\,\, \theta \in \Theta.
\end{align}
Equivalently, the subalgebra ${\cal \widetilde{A}}=S^{-1}(\cal B)$ of $\cal A$ is complete sufficient if any $\cal {\widetilde{A}}$-measurable
$\chi$ with $\displaystyle E_{\theta}\chi\equiv 0$ implies \eqref{eq:com}.
The concept of completeness is due to  Lehmann and {Scheff\'{e} and its role in the estimation theory is explained by the following.
%\begin{thm}
If $ \widetilde{\cal A}\,\,\,({\rm resp.},\: S(X))$ is a complete sufficient subalgebra (resp., statistic) for $\cal P$, then $\displaystyle \cal{\widetilde{A}}$-measurable statistics (resp., statistics depending on $X$ only through $S(X)$) with finite second moment and only they are UMVUEs.
%\end{thm}

Notice that the factorization theorem does not distinguish between complete and incomplete sufficiency. Proving (or disproving) completeness of a sufficient statistic requires an additional analysis. \cite{LS1950} (see also \cite{Lehmann1998}) showed that for the natural exponential families of full rank{\color{blue},} the minimal sufficient statistic is complete.

\subsection{A case of partial completeness}
\label{subsec:22}
The following interesting observation is due to \cite{B1983}. Let ${\cal{P}}=\{{P_{\theta,\eta}}\,,\,(\theta, \eta)\in  {\Theta \times \Xi}\}$ be a family of distributions on $(\cal X, \cal A)$ parameterized by a bivariate parameter. Suppose that for any fixed $\eta^{*} \in \Xi$, a subalgebra $\cal {\widetilde{A}}$ (resp., a statistic $S(X)$) is complete sufficient for $\theta$ (i.e., for the family ${\cal{P}^{*}}=\{P_{\theta,\eta^{*}},\,\theta\in \Theta\}$).
Then any $\displaystyle \cal{\widetilde{A}}$-measurable statistic (resp., depending on $X$ through $S(X)$) with finite variance is a UMVUE.
\noindent
The following models are covered by the Bondesson's result.
Suppose that the probability density function of $X=(X_1, X_2)$ is factorized as $$f\left(x_1,x_2;\, \theta, \eta\right)=R_1(x_1;\theta)R_2(x_2;\eta)r(x_1; x_2).$$
For any fixed $\eta=\eta^{*}$, $X_1$ is sufficient for $\theta$. If $X_1$ is complete sufficient, then any statistic $\widehat{g}(X_1)$ is a UMVUE. One can notice that,
in general, the parametric function $\displaystyle E_{\theta,\eta}\widehat{g}(X_1)=g(\theta, \eta)$ depends on the whole parameter $(\theta, \eta)$ and not only on the first component. The original example due to \cite{B1983} is a sample $\displaystyle X=\left(X_1,\ldots,X_n\right)$ from a gamma population with density function
$$\displaystyle f(x;\,\theta, \eta)=\frac{c(\alpha,\beta)}{\theta}(x-\eta)^{\alpha-1}\displaystyle e^{-\frac{x-\eta}{\beta \theta}},\,x>\eta$$
with $(\theta, \eta)\in \mathbb{R}_{+}\times \mathbb{R}$ as parameters, and $\alpha, \beta$ ($\alpha>0, \alpha\neq 1, \beta>0$) known.
For any fixed $\eta=\eta^{*},\,\overline{X}$ is complete sufficient for $\theta$. Any statistic $\widehat{g}(\overline{X})$ (with finite second moment) is the UMVUE for
$\displaystyle E_{\theta, \eta}\,\widehat{g}(\overline{X})=g(\theta, \eta)$. In particular, $\overline{X}$ is the UMVUE for $\displaystyle E_{\theta, \eta}(\overline{X})=\eta+\alpha \beta \theta.$

Actually, \cite{B1983} arguments prove a more general result.  Namely, if $\displaystyle \widehat{A}$ is an MVE-algebra for $\cal{P}^{*}$ (not necessarily a complete sufficient subalgebra) for any fixed $\eta^{*}$, then $\displaystyle \widehat{A}$
is an MVE-algebra for $\cal P$. Indeed, if $\displaystyle \chi\left(X\right)$ is a zero-mean statistic with finite second moment for $\cal P$, it is zero-mean statistic for $\cal P^{*}$. So that due to the assumption any $\cal{\widehat{A}}$-measurable $\widehat{g}(X)$ with $\displaystyle E_{\theta,\,\eta^{*}}|\widehat{g}(X)|^2<\infty$ one has
\begin{equation}
\label{eq:pc}
E_{\theta,\,\eta^{*}}\left(\widehat{g}\chi\right)=0,\,\,\,\,\forall\, \theta \in \Theta.
\end{equation}
Since \eqref{eq:pc} holds for any $\eta^{*}\in \Xi$,\,\, $\widehat{g}$ is a UMVUE for $\cal P$.
One can notice that the reverse is not true, in general. If $\widehat{g}$ is a UMVUE for $\cal P^{*}$, it is not necessarily a UMVUE for
$\cal P$.
The thing is that the relation
$\displaystyle E_{\theta,\,\eta^{*}}\chi(X)=0$ for some $\eta^{*}\in \Xi$ and
all $\theta \in \Theta$ does not imply $\displaystyle E_{\theta,\,\eta}\chi(X)=0$
for all $\eta\in \Xi,\,\,\,\theta \in \Theta$, so there are more zero mean statistics for $\cal P^{*}$ than for $\cal P$, in general.

\subsection{UMVUEs from categorical data}
\label{subsec:23}
Let $X$ be a categorical random variable whose values may be taken as $1,2,\ldots,N$. The distribution of $X$ is given by $$P\left(X=k;\theta\right)=p_{k}(\theta),\,\,k=1,\ldots, N$$ with $\theta \in \Theta$ as a parameter.
Plainly, only parametric functions from $\displaystyle L=\text{span}\left\{p_1(\theta), p_2(\theta),\ldots, p_N(\theta)\right\}$ are estimable, i.e.,
can be unbiasedly estimated. The set $\displaystyle M=\left\{p_1(\theta), p_2(\theta),\ldots, p_N(\theta)\right\}$ can be partitioned into
$$M=M_1 \cup \ldots \cup M_r$$ in such a way that
\begin{itemize}
\item[(i)]
the subspaces $L_1=\text{span}\left\{M_1\right\}, \ldots, L_r=\text{span}\left\{M_r\right\}$ are linearly independent, and
\item[(ii)]
the partition is maximal, i.e., if for some $l$, $M_l=M_l^{'}\cup M_l^{''}$, $M_l^{'}\cap M_l^{''}=\emptyset$ and (i) holds for the new partition, then either $M_l=M_l^{'}$ or $M_l=M_l^{''}$.
\end{itemize}
Such a partition is unique up to ordering. Without loss of generality, one may assume
\begin{align}
\label{eq:111}
&
M_1=\left\{p_1(\theta),\ldots, p_{k_1}(\theta)\right\},\,\,\,\,M_2=\left\{p_{k_1+1}(\theta),\ldots, p_{k_1+k_2}(\theta)\right\},\ldots,\nonumber\\
&
 M_r=\left\{p_{k_1+\ldots+k_{r-1}+1}(\theta),\ldots, p_{N}(\theta)\right\}.
\end{align}

The partition \eqref{eq:111} generates a partition of the set $\displaystyle \left\{1,2,\ldots,N\right\}$:
\begin{align}
\label{eq:par}
&
I_{1}\cup I_{2} \cup \ldots \cup I_{r}=
\left\{1,\ldots,k_1\right\}\cup \left\{k_1+1,\ldots,k_1+k_2\right\}\cup \ldots \cup \left\{k_1+\ldots+k_{r-1}+1,\ldots,N\right\}.
\end{align}
A statistic $\displaystyle \widehat{g}(X)$ is a UMVUE if and only if it is constant on the elements of \eqref{eq:par},
$$\widehat{g}(x)=const=g_j,\,\,\,x\in I_j,\,\,j=1,\ldots,r.$$
The subalgebra of UMVUEs is generated by sets \eqref{eq:par}. The class of parametric functions admitting UMVUEs is a linear subspace of $L$,
\begin{align}
&
\text{span}\,\left\{\pi_1(\theta),\ldots,\pi_r(\theta)\right\},\,\,\,\pi_j(\theta)=\sum_{k \in I_{j}}p_{k}(\theta),\,\,\,\,j=1,\ldots,r.
\end{align}
Note that $S(X)$ is the minimal sufficient statistic for $\theta$ if and only if $\displaystyle S(k)=S(l)$ is equivalent to
\begin{align}
\label{eq:ms}
&
p_k(\theta)=c_{kl}p_l(\theta),\,\,\,\,\,\forall \theta \in \Theta
\end{align}
for some constant $c_{kl}>0$.
Due to (i), a pair $(k,l)\in \left\{1,2,\ldots,N\right\}$ with \eqref{eq:ms} always belongs to the same element of the partition \eqref{eq:par}.
The above construction is due to \cite{KK2006}. The following example with
$$p_1(\theta)=\theta,\,\,p_2(\theta)=\theta^2,\,\,p_3(\theta)=\theta+\theta^2,\,\,p_4(\theta)=1-2\theta-2\theta^2,\,\,\,\theta\in\left(0,1/4\right)$$
illustrates  the situation when the minimal sufficient statistic is trivial and incomplete while the subalgebra of UMVUEs is generated by two sets, $\left\{1,2,3\right\}$
and $\left\{4\right\}$. Therefore only elements of $\text{span}\left\{1, \theta+\theta^2\right\}$ possess UMVUEs.

\subsection{Lehmann's example}
\label{subsec:24}
The following example is due to \cite{LS1950} (see also \cite{Lehmann1998}, pp. 84--85).
Let $X \in \{-1,0,1,2,\ldots\}$ with
$$P_{\theta}\left(X=-1\right)=\theta,\,\,\,\, P_{\theta}\left(X=k\right)=(1-\theta)^2\theta^k,\,\,\,\theta \in (0,\,1),\,\, k=0,1,2,\ldots$$
It is easy to see that all unbiased estimators of zero are of the form
$$U(X)=aX\,\,\,\,\text{for some}\,\,\,\,a\in \mathbb{R}.$$
If $\widehat{g}(X)$ is a UMVUE, then
$$E_{\theta}\left\{\widehat{g}(X)X\right\}=0,\,\,\,\theta\in(0,1)$$
so that $X\widehat{g}(X)$ is itself an unbiased estimator of zero. Thus,
$$x\widehat{g}(x)=ax,\,\,\,x\in\{-1,0,1,2,\ldots\}\,\,\,\,\text{for some}\,\,\,\,a\in \mathbb{R},$$
whence $$\widehat{g}(x)=\widehat{g}(-1)\,\,\,\text{for all}\,\,\,x\neq 0,\,\,\,\text{with an arbitrary}\,\,\,g(0).$$
The subalgebra of UMVUEs in this example is generated by two sets, $\{0\}$ and  $\{-1,1,2,\ldots\}$. The parametric functions possessing UMVUEs are elements of $\text{span}\{1, (1-\theta)^2\}$.
\section{Properties of MVE-algebras}
\label{sec:3}
In this section some properties of MVE-algebras are presented. They are similar to properties of complete sufficient statistics/subalgebras and, in our opinion, are an argument in favor of that under rather general conditions on the statistical model $\displaystyle \left(\mathcal{X}, \mathcal{A}, \mathcal{P}\right)$, there exists a maximal MVE-algebra $\widehat{\cal{A}}$ such that a statistic $\widehat{g}(X)$ with $E_{\theta}|\widehat{g}(X)|^{2}<\infty,\,\,\forall\,\,\theta \in \Theta$ is a UMVUE if and only if it is
$\widehat{\cal{A}}$-measurable. \\
\noindent
Recall that a subalgebra $\widehat{\cal{A}}\subset \cal{A}$ is an MVE-algebra  if any $\widehat{\cal{A}}$-measurable statistic $\widehat{g}(X)$ with finite second moment is a UMVUE.
\begin{thm}
If $\widehat{\cal{{A}}}_{1},\widehat{\cal{A}}_{2}$ are MVE-algebras, so is $\widehat{\cal{A}}=\sigma(\widehat{\cal{A}}_{1},\:\widehat{\cal{A}}_{2})$, the smallest $\sigma$-algebra containing both $\widehat{\cal{A}}_{1}$ and $\widehat{\cal{A}}_{2}$ (in other words, generated by  $\widehat{\cal{A}}_{1}$ and $\widehat{\cal{A}}_{2}$).
\end{thm}
In terms of statistics, if $U_{1}(X),\:U_{2}(X)$ are such that any estimators $\hat g_{1}(U_{1}(X)),\:
\hat g_{2}(U_{2}(X))$ with finite second moments are UMVUEs, so is $\hat g(U_{1}(X),\:U_{2}(X))$ with finite second moment.
\begin{proof}
For $A_1 \in \widehat{\cal{{A}}}_{1}, A_2 \in \widehat{\cal{A}}_{2}$ set $h_{i}(x)={1}_{\left\{x\in A_i\right\}}, i=1,2.$ Since $h_1(X)$ is a UMVUE,
${\rm cov}_{\theta}\left(h_1(X), \chi(X)\right)=E_{\theta}\left(h_1(X), \chi(X)\right)=0, \forall\,\,\theta \in \Theta$ for any $\chi(X)$ with
$E_{\theta}\left(\chi(X)\right)=0, E_{\theta}|\chi(X)|^2<\infty\,\, \forall\,\,\theta \in \Theta$. Thus $\chi_{1}(X)=h_1(X)\chi(X)$ is also a zero mean statistic. Since
$|h_{1}(X)|\leq 1,\,E_{\theta}|\chi(X)|^2<\infty$,
$${\rm cov}_{\theta}\left(h_2(X), \chi_1(X)\right)=E_{\theta}\left(h_2(X)h_1(X)\chi(X)\right)=0,\,\,\, \forall\,\,\theta \in \Theta$$ and
$h_1(X)h_2(X)=1_{\left\{X\in A_1\cap A_2\right\}}$ is a UMVUE.\\
Let now $A_{1i} \in \widehat{\cal{{A}}}_{1}, A_{2i} \in \widehat{\cal{A}}_{2}, i=1,\ldots,n$ with indicators $h_{1i}(x)={1}_{\left\{x\in A_{1i}\right\}}, h_{2i}(x)={1}_{\left\{x\in A_{2i}\right\}}, i=1,\ldots,n$. The above arguments prove that for any constants $c_1,\ldots,c_n$
\begin{equation}
\label{eq:U}
\sum_{i=1}^{n}c_ih_{1i}(X)h_{2i}(X)
\end{equation}
is a UMVUE.  As is well known, the functions \eqref{eq:U} are dense in the Hilbert space $L_{\theta}^{2}\left(\widehat{\cal{A}}\right)$ of
$\widehat{\cal{A}}$-measurable functions $h(X)$ with $\int |h(x)|^2 dP_{\theta}(x)<\infty.$
For an $\widehat{\cal{A}}$-measurable statistic $h(X)$ with finite second moment and given $\varepsilon>0$  take  $\widehat{h}(X)$ of the form \eqref{eq:U}
with $E_{\theta}|h(X)-\widehat{h}(X)|^{2}\leq \varepsilon^2$. Now, for any unbiased estimator of zero $\chi(X)$ with finite second moment,
\begin{align*}
|{\rm cov}_{\theta}\left(h(X),\chi(X)\right)|=|{\rm cov}_{\theta}\left(\widehat{h}(X),\chi(X)\right)+{\rm cov}_{\theta}\left(h(X)-\widehat{h}(X),\chi(X)\right)|
\leq \varepsilon \sqrt{E_{\theta}|\chi(X)|^2}.
\end{align*}
Since, $\varepsilon>0$ is arbitrary, ${\rm cov}_{\theta}\left(h(X),\chi(X)\right)=0$ and $h(X)$ is a UMVUE.
\end{proof}
 The next claim is an immediate corollary of Theorem 2.
\begin{cor}
 Let $\left\{\widehat{\cal{A}}_{\gamma}, \gamma \in \Gamma \right\}$ be the collection of all MVE-algebras indexed by $\gamma \in \Gamma$. Then $\widehat{\cal{A}}=\sigma(\widehat{\cal{A}}_{\gamma},\:\gamma \in \Gamma)$ is the maximal MVE-algebra.
\end{cor}
Turn now to the setup of combining independent data. Let $({\cal X}_i,\:{\cal A}_i,\:{\cal P}_i)$
be statistical models with ${\cal P}_i =\{P_{\theta_i},\:\theta_i \in \Theta_i\},\:i=1,\:2$.
Set
\[({\cal X},\:{\cal A},\:{\cal P})=({\cal X}_1 \times {\cal X}_2,\:{\cal A}_1 \otimes {\cal A}_2,\:
{\cal P}_1 \times {\cal P}_2)\]
is a statistical model with ${\cal P}=\{P_{\theta}=P_{\theta_1} \times P_{\theta_2}\}$ parameterized by a ``bivariate'' parameter $\theta =(\theta_1,\:\theta_2)\in \Theta_1 \times\Theta_2 =\Theta$.
\begin{thm}
If $\widehat {\cal A}_i \subset {\cal A}_i$ is an MVE-algebra in the model $({\cal X}_i,\:{\cal A}_i,\:{\cal P}_i),\:i=1,\:2$, then $\widehat {\cal A}=\widehat{\cal A}_1 \otimes \widehat{\cal A}_2$ is an MVE-algebra in the model $({\cal X},\:{\cal A},\:{\cal P})$.
\end{thm}
In terms of statistics, let $X_i \sim P_{\theta_i},\:i=1,\:2$ and $X=(X_1,\:X_2)\sim P_{\theta}$.
If statistics $U_{i}(X_i),\:i=1,\:2$ are such that any estimator $\hat g_{i}(U_{i}(X_i))$ with a finite second moment is a UMVUE of a parametric function $g_{i}(\theta_i), i=1,\:2$, then
any estimator $\hat g(U_{1}(X_1),\:U_{2}(X_2))$ with finite second moment is a UMVUE of some $g(\theta_1,\theta_2).$
\begin{proof} In view of Theorem 2, suffice to prove that $\widehat{\cal A}_1$ is an MVE-algebra in the model  $({\cal X},\:{\cal A},\:{\cal P})$. Let $\chi(X_1,\:X_2)$ be a zero mean statistic with finite second moment so that
\[\int_{{\cal X}_{1}} \int_{{\cal X}_{2}} |\chi(x_1,\:x_2)|dP_{\theta_1} (x_1)dP_{\theta_2} (x_2)<\infty, \:(\theta_1,\:\theta_2)\in\Theta.\]
By Fubini theorem, for any $\theta_2$ the function
\[\tilde\chi(x_1;\:\theta_2)=\int_{{\cal X}_{2}} \chi(x_1,\:x_2) dP_{\theta_2}(x_2)\]
is well defined for $P_{\theta_1}$-almost all $x_1$ and is a zero mean statistic for the model
$({\cal X}_1,\:{\cal A}_1,\:{\cal P}_1)$. Thus, for any $\widehat{\cal A}_1$-measurable statistic
$U_{1}(X_1)$ with finite second moment,
\[0={\rm cov}(U_{1}(X_1),\:\tilde\chi(X_1;\:\theta_2))=\int_{{\cal X}_{1}} U_{1}(x_1)\tilde\chi(x_1;\theta_2) dP_{\theta_1}(x_1)=\int_{{\cal X}_{1}} \int_{{\cal X}_{2}}
 U_{1}(x_1) \chi(x_1,\:x_2)dP_{\theta_1} (x_1)dP_{\theta_2}\]
and $U_{1}(X_1)$ is a UMVUE in the combined model  $({\cal X},\:{\cal A},\:{\cal P})$.
\end{proof}
If in Theorem 3{\color{blue},} the condition ``$\widehat{\cal A}_i \subset {\cal A}_i$ is an MVE-algebra for $\theta_i,\:i=1,\:2$''
is replaced with ``$\widehat{\cal A}_i \subset {\cal A}_i$ is a complete sufficient subalgebra'',
then the claim may be replaced with ``$\widehat {\cal A}=\widehat{\cal A}_1 \otimes \widehat{\cal A}_2$ is a complete sufficient subalgebra for $(\theta_1,\:\theta_2)$''. This result was proved
in \cite{LR1976} (where it is stated for arbitrary subalgebras, not necessarily sufficient) strengthening previous results by \cite{Plachky1977} and \cite{Fraser1957} that required some additional conditions on the models.\\

\section{Comments}
In conclusion we would like to make some general comments on sufficiency and complete sufficiency.
While sufficiency plays a fundamental role in all areas of statistical inference, complete sufficiency is tailored for the estimation. In a sense, while sufficiency is a statistical
concept, completeness looks more like a mathematical-statistical tool.
It seems that a more general concept of an MVE-algebra preserves the basic property of complete sufficiency in some setups when the (minimal) sufficient statistic is incomplete.

%%%%%%%%%%%%%%%%%%%%%%%%%%%%%%%%%%%%SUPPLEMENT%%%%%%%%%%%%%%%%%%%%%%%%%%%%%%%5%%%%%
\end{document}